\documentclass{amsart}

\usepackage[utf8]{inputenc}
\usepackage[T1]{fontenc}
\usepackage{graphicx}
\usepackage{subfigure}
\usepackage{epsfig}
\usepackage{geometry}
\usepackage{float}
\usepackage{empheq,etoolbox}
\usepackage[num,bibjustif,recuo=0.8cm]{abntex2cite} 
\citebrackets[] % coloca os números das citações entre [ ]

\usepackage{amsthm}

\usepackage{amsmath,amssymb,exscale}

\usepackage{framed}

\usepackage{enumerate}

\newtheorem{teo}{Theorem}[section]
\newtheorem*{teoa}{Theorem A}
\newtheorem*{teob}{Theorem B}

\newtheorem{ex}[teo]{Example}

\newtheorem{prop}[teo]{Proposition}
\newtheorem{coro}[teo]{Corollary}
\newtheorem{lema}[teo]{Lemma}

\newcommand{\dem}{\noindent \underline{Proof.}: $\,$}
\newcommand{\fim}{\hfill $\rule{2.0mm}{2.0mm}$ \\}    % marca��o de fim das demonstra��es

\newcommand{\h}{\mathcal{H}}
\newcommand{\G}{\mathcal{G}}
\newcommand{\F}{\mathcal{F}}
\newcommand{\W}{\mathcal{W}}
\newcommand{\leg}{\operatorname{Leg}}
\newcommand{\tang}{\operatorname{Tang}}
\newcommand{\sing}{\operatorname{Sing}}

\title{Webs Generated by Products of convex and homogeneous Foliations on $\mathbb{P}^2$}
\author{Carla Pracias \and Maycol Falla Luza} 
\email{carla\_pracias@id.uff.br \and hfalla@id.uff.br}
\date{}
\begin{document}
\begin{abstract}
This paper investigates flat webs on the projective plane. We present two methods for constructing such webs: the first involves taking the product of finitely many convex reduced foliations and invariant lines, while the second consists of taking the product of finitely many convex homogeneous foliations and invariant lines. In both cases, we demonstrate that the dual web is flat.
\end{abstract}

%==========================================================================
% Início do texto do artigo
%==========================================================================

\maketitle

\setcounter{secnumdepth}{-1} %para capitulos sem numeração

\nocite{S2} \nocite{S} \nocite{B} \nocite{M} \nocite{ordinary} \nocite{JD} \nocite{J} %ordem alfabética da referência

\section{Introduction}

Locally, in affine coordinates on the projective plane, a foliation is defined by a polynomial $1$-form $\omega=a(x,y)dx+b(x,y)dy$ with isolated zeros. Similarly, a $k$-web is defined by a $k$-symmetric polynomial $1$-form with non identically zero discriminant and isolated zeros as follows
$$\Omega = \sum\limits_{i+j=k}a_{ij}(x,y)dx^idy^j.$$

In the context of this article, we are interested in the curvature of webs. Given a web $\W$ on a complex surface, its curvature is a meromorphic $2$-form with poles contained in the discriminant $\Delta(\W)$, denoted by $K(\W)$. This $2$-form satisfies $\varphi^*K(\W)=K(\varphi^*\W)$ for any holomorphism $\varphi$. When the curvature is zero we say that $\W$ is flat. A key observation: since there are no holomorphic $2$-forms on $\mathbb{P}^2$, the web is flat if and only if the curvature is holomorphic at the generic points of the irreducible components of $\Delta(\W)$. A important result of web geometry, due to Blaschke-Dubourdien, characterizes the local equivalence of a (germ of) $3$-web $\W$ on $\mathbb{C}^2$ with the trivial $3$-web defined by $dx\cdot dy\cdot(dx-dy)$ through the vanishing of curvature of $\W$. The usually definition of curvature for a $k$-web $\W$ with $k>3$  is the sum of the curvatures of all $3$-subwebs of $\W$, and we say that the $k$-web is flat when this sum is zero. To our knowledge there is not a characterization of the flatness of a web,  nevertheless is a necessary condition for the maximality of the rank of the web, according to a result of Mihaileanu.\\

For a global context, given a $k$-web $\W$ of degree $d$ on the projective plane we can consider the Legendre transform of $\W$, it will be denoted by $\leg\W$. This is a $d$-web of degree $k$ on the dual projective plane $\Check{\mathbb{P}}^2$ which we explain now. Take a generic line $l$ on $\mathbb{P}^2$ and consider the tangency locus $\tang(\W,l)=\{p_1,\hdots,p_d\}\subset \mathbb{P}^2$ between $\W$ and $l$. We can think the dual $\Check{p}_i$ as lines on $\Check{\mathbb{P}}^2$ passing through the point $l$ of $\Check{\mathbb{P}}^2$. Then the set of tangent lines of $\leg\W$ at $l$ is just $T_l\leg\W=\displaystyle\bigcup_{i=1}^d\Check{p}_i$. More precisely, let $(p,q)$ be the affine chart of $\Check{\mathbb{P}}^2$ correspond to the line $\{y=px+q\}\subset\mathbb{P}^2$. If $\W$ is defined by an implicit affine equation $F(x,y;p)=0$ with $p=dy/dx$ then $\leg \W$ is given by the implicit differential equation 
$$\Check{F}(p,q;x):=F(x,px+q;p), \hspace{1cm} \text{with } \hspace{1cm} x=-\frac{dq}{dp}.$$

 When the web $\W$ decomposes as the product of two webs $\W_1\boxtimes\W_2$ then we have that the following decomposition also occurs: $\leg\W=\leg\W_1\boxtimes\leg\W_2$. For more details about the definition of Legendre transform, see for instance \cite[Part I Chapter II Section 2.5]{ordinary}.\\

In particular, we are especially interested in convex foliations. For us a foliation $\F$ on $\mathbb{P}^2$ is convex if its inflection divisor is completely invariant by $\F$. Even more, if this divisor is also reduced we will say that $\F$ is reduced convex foliation. \\

Our first result is about the dual of a web generated by product of some invariant lines and  reduced convex foliations under some hypotheses. 

\begin{teoa}\label{teoA}
    Let $\F_1, \hdots, \F_n$ be reduced convex foliations on $\mathbb{P}^{2}$ of degrees $d_1, \hdots, d_n$ such that $\tang(\F_i,\F_j)$ is formed by invariant lines by $\F_i$ and $\F_j$ where $i\neq j$ in $\{1,\hdots,n\}$. If $l_1, \hdots, l_k$ are invariant lines by all the foliations then the transform $\leg \W$ of the web $\W=l_1\boxtimes\cdots\boxtimes l_k\boxtimes\F_1\boxtimes\cdots\boxtimes\F_n$ is flat.
\end{teoa}

We present some reduced convex foliations that satisfies the theorem conditions and provide a family of examples of flat webs of this type. 

In the second part of this paper, we will work with homogeneous foliations. A foliation $\h$ of degree $d$ on $\mathbb{P}^2$ is called homogeneous if there exists a system of affine coordinates $(x,y)$ where $\h$ is defined by a homogeneous vector field $A(x,y)\frac{\partial}{\partial x}+B(x,y)\frac{\partial}{\partial y}$ with $A,B\in \mathbb{C}[x,y]_d$ and $gcd(A,B)=1$. In the same direction as the first result, we want to determine the flatness of the Legendre transform of a web generated by the product of invariant lines and many convex homogeneous foliations. The result can
be phrased as follows.

\begin{teob}\label{teoB}
    Let $\h_1, \hdots, \h_n$ be convex homogeneous foliations on $\mathbb{P}^{2}$ of respective degrees $d_i$ with $i=1,\hdots, n$, such that, $\tang(\h_i,\h_j)=L_\infty \cup C_{ij}^{inv}$ where $\forall i\neq j$ in $\{1,\hdots,n\}$, $C_{ij}^{inv}$ is a product of invariant lines by $\h_i$, $\h_j$. If $l_1, \hdots, l_k$ are invariant lines by all the foliations then the dual web $\leg \W$ of $\W=l_1\boxtimes\cdots\boxtimes l_k\boxtimes\h_1\boxtimes\cdots\boxtimes\h_n$ is flat.  
\end{teob}

These results are expansion of those presented by Bedrouni in \cite{S2}. The Theorem A can be seen as a generalization of \cite[Theorem 1]{S2}, while the Theorem B is a generalization of \cite[Theorem 2]{S2}. Furthermore, note that the hypothesis that the tangency locus of two foliations be formed by invariant lines is necessary, as can be seen in the Example \ref{ex3}.

%%%%%%%%%%%%%%%%%%%%%%%%%%%%%%%%%%%%%%%%%%%%%%%%%%%%%%%%%%%
\setcounter{secnumdepth}{1}
\section{Preliminaries}
Let $\F$ be a satured foliation on $\mathbb{P}^2$, we can define the Gauss map of $\F$ as the rational map $\G_\F:\mathbb{P}^2 \dashrightarrow\Check{\mathbb{P}}^2$ such that $\G_\F(p)=T_p\F$, which is well defined outside the singular set of $\F$, $\sing\F$. If the homogeneous form $\omega=a(x,y,z)dx+b(x,y,z)dy+c(x,y,z)dz$ defines the foliation $\F$, then the Gauss map can be described as
$$\G_\F(p)=[a(p):b(p):c(p)].$$

The inflection divisor of $\F$, denoted by $I(\F)$, is defined by the vanishing of the determinant
$$\det \left (\begin{matrix}
    x && y && z\\
    X(x) && X(y) && X(z)\\
    X^2(x) && X^2(y) && X^2(z)
\end{matrix} \right )$$
with $X$ being the vector field defining $\F$. This divisor has some relevant properties for us:
\begin{enumerate}[(i)]
    \item $I(\F)$ does not depend on a particular choice of system of homogeneous coordinates on $\mathbb{P}^2$.
    \item On $ \mathbb{P}^2 \setminus \sing(\F)$, $I(\F)$ coincides with the curve described by the inflection points of the leaves of $\F$.
    \item If $C$ is an irreducible algebraic invariant curve of $\F$ then $C\subset I(\F)$, if and only if, $C$ is an invariant line.
    \item The degree of $I(\F)$ is exactly $3d$, where $d$ is the degree of $\F$.
\end{enumerate}

For a discussion in a more general context and additional details about inflection divisor, see \cite{J}. The following result will be useful for us.

\begin{lema}[\cite{M}, Lemma 2.2]
    Let $\F$ be a foliation on $\mathbb{P}^2$. Then
    $$\Delta(\leg\F)=\G_\F(I(\F))\cup\Check{\Sigma}(\F)$$
    where $\Check{\Sigma}(\F)$ consists in the dual lines of the special singularities 
    $$\Sigma(\F)=\{s\in\sing(\F):\nu(\F,s)\ge2 \text{ or } s \text{ is a radial singularity of } \F\}$$
    of $\F$ and $\nu(\F,s)$ stands for the algebraic multiplicity of $\F$ at $s$.
\end{lema}

We also denote by $\Sigma_{\F}^{\text{rad}}$ the set of radial singularities and by $\Check{\Sigma}_{\F}^{\text{rad}}$ the dual lines. It is convenient to define a similar nomenclature: we will denote by $\Sigma_{\F}^{l}$ where $l$ is a line, the set of singularities of $\F$ in $l$ i.e. $\sing\F\cap l$ and by $\Check{\Sigma}_{\F}^{l}$ the respective dual lines. In order to study the product of foliations, we shall need the following proposition.

\begin{prop}\label{teo1} 
Let $\W=\F_1 \boxtimes \F_2$ be a $(d_1+d_2)$-web in $\mathbb{P}^2_{\mathbb{C}}$, where $\F_1$ and $\F_2$ are foliations of respective degrees $d_1$ and $d_2$. If $\tang(\F_1,\F_2)=C^{\text{inv}}\cup C^{\text{tr}}$, where $C^{\text{inv}}$ is the invariant part for both $\F_1$ and $\F_2$, and $C^{\text{tr}}$ is the transversal part, then

$$\Delta(\leg\W)=\Delta(\leg \F_1)\cup\Delta(\leg \F_2)\cup\G_{\F_i}(C^{\text{inv}})\cup \G_{\F_i}(C^{\text{tr}})$$
for each $i=1,2$.
\end{prop}

\dem We know that $\Delta(\leg \W)=\Delta(\leg \F_1)\cup\Delta(\leg \F_2)\cup \tang (\leg \F_1,\leg \F_2).$ Clearly we have $\G_{\F_i}(C^{\text{inv}})\cup \G_{\F_i}(C^{\text{tr}}) \subset \tang (\leg \F_1,\leg \F_2)$. Now, take an irreducible curve $C \subseteq \tang(\leg \F_1,\leg \F_2)\setminus (\Delta(\leg\F_1)\cup\Delta(\leg\F_2))$ and $\Check{l} \in C$ a generic point. Then there is a point $p\in l$ such that $p \in \tang(\F_i,l)$ for $i=1,2$ and we would have that $p \in \tang(\F_1,\F_2)$ and $\G_{\F_i}(p)=\Check{l}$. Thus, $C \subseteq \G_{\F_i}(C^{\text{inv}})\cup \G_{F_i}(C^{\text{tr}})$. \fim

For our results, we will need the following lemma.

\begin{lema}\label{lema1}
Let $\W_1,\W_2, \ldots \W_k$ be $\nu_i$-webs and let $C \subset \Delta(\W_i)$ be an irreducible component totally invariant by $\W_i$ and having minimal multiplicity $\nu_i - 1$ for $i=1,\ldots, k$. Let $\W_d$ be a smooth $d$-web transverse to $C$. Then the curvature of $\W= \W_1 \boxtimes \ldots \boxtimes \W_k \boxtimes \W_d$ is holomorphic along $C$.
\end{lema}

\dem 
The proof is a simple recurrence of the case $k=2$. So let us assume we are in this case.

Suppose that $1 \le \nu_1 \leq \nu_2$. Let $(U,(z,w))$ be a local coordinate system such that $C\cap U=\{w=0\}$ and 
$$T\W_1|_{U}=\{dw^{\nu_1}+w(a_{\nu_1-1}(z,w)dw^{\nu_1-1}dz + \hdots +a_0(z,w)dz^{\nu_1})=0\} \hspace{1cm} \text{with } w\nmid a_0$$

$$T\W_2|_{U}=\{dw^{\nu_2}+w(b_{\nu_2-1}(z,w)dw^{\nu_2-1}dz + \hdots +b_0(z,w)dz^{\nu_2})=0\} \hspace{1cm} \text{with } w\nmid b_0$$

and

$$T\W_d|_{U}=\{ \prod_{j=1}^{d}(dz + \frac{1}{\nu}g_j(z,w)dw)\}$$

where $\nu=\nu_1\cdot\nu_2$. Let $\pi_i: \Bar{U} \rightarrow U$ be the ramified covering given by $(z,w)=\pi(x,y)=(x,y^{\nu_i})$ for $i=1,2$ and $\pi=\pi_1\circ\pi_2=\pi_2\circ\pi_1=(x,y^\nu)$.

The irreducibility of $\W_1$ and $\W_2$ implies that its monodromy group is cyclic and consequently $\pi_i^{*}\W_i$ is totally decomposable. More precisely, $\pi_i^{*}\W_i$ is given by

$$y^{\nu_1(\nu_1-2)}dy^{\nu_1} +\Bar{a}_{\nu_1-1}(x,y)y^{(\nu_1-1)^2}dy^{\nu_1-1}dx+\hdots +\Bar{a}_0dx^{\nu_1}=0$$

$$y^{\nu_2(\nu_2-2)}dy^{\nu_2} +\Bar{b}_{\nu_2-1}(x,y)y^{(\nu_2-1)^2}dy^{\nu_2-1}dx+\hdots +\Bar{b}_0dx^{\nu_2}=0$$

Since $y\nmid\Bar{a}_0$ and $y\nmid\Bar{b}_0$ we can write

$$\pi_i^{*}\W_i : \prod_{j=1}^{\nu_i}(dx + y^{\nu_i-2}f_i(x,\xi_i^j y)\xi_i^{-j}dy)=0$$

where $\xi_i=e^{\frac{2\pi i}{\nu_i}}$ and $y \nmid f_i$ for each $i=1,2$. Now, we have

$$\pi^{*}\W_1 = \pi_2^{*}(\pi_1^{*}\W_1) :  \prod_{j=1}^{\nu_1}(dx + \nu_2 y^{\nu-\nu_2-1}f_1(x,\xi_1^j y^{\nu2})\xi_1^{-j}dy)=0$$

$$\pi^{*}\W_2 = \pi_1^{*}(\pi_2^{*}\W_2) :  \prod_{j=1}^{\nu_2}(dx + \nu_1 y^{\nu-\nu_1-1}f_2(x,\xi_2^j y^{\nu1})\xi_2^{-j}dy)=0$$

and 

$$\pi^{*}\W_d : \prod_{j=1}^{d}(dx + y^{\nu-1}g_j(x,y^\nu)dy)=0$$

On the other hand, $K(\pi^{*}\W)=\sum\limits_{i<j<k}d\eta_{ijk}$ where $\eta_{ijk}$ is the unique $1$-form that satisfies $d(\delta_{rs}\omega_t)=\eta_{ijk}\wedge\delta_{rs}\omega_t$ for each permutation on $(r,s,t)$ of $(i,j,k)$ and the function $\delta_{rs}$ is defined by $\omega_r\wedge\omega_s=\delta_{rs}(x,y)dx\wedge dy$.

We denote by $\varphi(x,y)=(x,\xi y)$ where $\xi=e^\frac{2\pi i}{\nu}$ the deck transformation of $\pi$. We have

$$K(\pi^*\W)=\pi^*K(\W)=\varphi^*\pi^*K(\W)=\varphi^*K(\pi^*\W).$$

Note that, $\varphi^*(y^ndx\wedge dy)=y^ndx\wedge dy$ if and only if $n  \equiv -1  \mod \nu$. 

If $\sum\eta_{ijk}=Adx+Bdy$ and since $K(\pi^*(\W_2\boxtimes\W_d))$ is holomorphic by \cite[Proposition 2.6]{JD}, we conclude using \cite[Lemma 2.2]{JD}:

\begin{enumerate}[i.]
    \item The order of the poles of $A$ along $\{y=0\}$ is $\le \nu-\nu_2$;
    \item $B$ is logarithmic along $\{y=0\}$ with constant residue.
\end{enumerate}

Thus we can write $A=\sum\limits_{n=\nu_2-\nu}^{\infty}\alpha_n(x)y^n$ and $B=\sum\limits_{n=-1}^{\infty}\beta_n(x)y^n$,  where $b_{-1}$ is constant. Then $d(\sum\eta_{ijk})=(B_x-A_y)dx\wedge dy$ with polar part given by
$$\Omega=\sum\limits_{\nu_2-\nu-1}^{-2}(m+1)a_{m+1}y^m, \hspace{1cm} m=n-1$$

This fact is summarized with remarks (i) and (ii) which will imply that $d\eta_{ijk}$ is holomorphic along $\{y=0\}$ and consequently $K(\W)$ is holomorphic along $C$. 
\fim 

An adaptation of \cite[Proposition 2.9]{S} allows us to show a similar result with an additional tangent foliation.

\begin{lema}\label{lema2}
    Let $\W_1,\ldots, \W_k$ be $\nu_i$-webs and let $C \subset \Delta(\W_i)$ for $i=1,\ldots, k$ as in the previous lemma. Let $\F$ be a foliation having $C$ as an invariant curve and let $\W_d$ be a smooth $d$-web transverse to $C$. Then the curvature of $\W=\F\boxtimes\W_1\boxtimes \ldots \boxtimes\W_k\boxtimes\W_d$ is holomorphic along $C$.
\end{lema}

\subsection{Some general results about curvature}  In order to establish our main results, we shall need the following propositions about some decomposable webs.

\begin{prop}\label{prop4}
    Let $\W=\W_1\boxtimes\cdots\boxtimes\W_n$ be a product of webs. If for each $i,j,l=1,\hdots,n$ we have:
    \begin{enumerate}[i.]
        \item $K(\W_i)=0$;
        \item $K(\W_i\boxtimes\W_j)=0$, $i\ne j$;
        \item $K(\W_i\boxtimes\W_j\boxtimes\W_l)=0$ with $i$, $j$, $l$ pairwise distinct. 
    \end{enumerate}
    Then $K(\W)=0$.
\end{prop}

\dem
Firstly, we can decompose each web as $\W_i=\F_1^i\boxtimes\cdots\boxtimes\F_{n_i}^i$ in a neighborhood of a generic point outside the discriminant. By the hypothesis we have:
\begin{eqnarray*}
    0 = K(\W_i\boxtimes\W_j) &=& K(\W_i)+K(\W_j) + \sum\limits_{\overset{1\le k < r \le n_i}{1 \le s \le n_j}} K(\F_k^i\boxtimes\F_r^i\boxtimes\F_s^j) + \sum\limits_{\overset{        1\le k < r \le n_j}{1 \le s \le n_i}} K(\F_k^j\boxtimes\F_r^j\boxtimes\F_s^i)\\
    &=&\sum\limits_{\overset{1\le k < r \le n_i}{1 \le s \le n_j}} K(\F_k^i\boxtimes\F_r^i\boxtimes\F_s^j) + \sum\limits_{\overset{        1\le k < r \le n_j}{1 \le s \le n_i}} K(\F_k^j\boxtimes\F_r^j\boxtimes\F_s^i)
\end{eqnarray*}

and

\begin{eqnarray*}
    0 = K(\W_i\boxtimes\W_j\boxtimes\W_k) &=&  \sum\limits_{\overset{1\le k < r \le n_i}{1 \le s \le n_j}} K(\F_k^i\boxtimes\F_r^i\boxtimes\F_s^j) + \sum\limits_{\overset{        1\le k < r \le n_j}{1 \le s \le n_i}} K(\F_k^j\boxtimes\F_r^j\boxtimes\F_s^i) \\
    && + \sum\limits_{\overset{1\le k < r \le n_i}{1 \le s \le n_l}} K(\F_k^i\boxtimes\F_r^i\boxtimes\F_s^l) + \sum\limits_{\overset{1\le k < r \le n_l}{1 \le s \le n_i}} K(\F_k^l\boxtimes\F_r^l\boxtimes\F_s^i)\\
    && + \sum\limits_{\overset{        1\le k < r \le n_l}{1 \le s \le n_j}} K(\F_k^l\boxtimes\F_r^l\boxtimes\F_s^j) + \sum\limits_{\overset{1\le k < r \le n_j}{1 \le s \le n_l}} K(\F_k^j\boxtimes\F_r^j\boxtimes\F_s^l)\\
    &&+ \sum\limits_{i<j<l}\sum\limits_{k,r,s=1}^{n_i,n_j,n_k}K(\F_k^i\boxtimes\F_r^j\boxtimes\F_s^l) + K(\W_i)+K(\W_j) + K(\W_k)\\    &=&\sum\limits_{i<j<l}\sum\limits_{k,r,s=1}^{n_i,n_j,n_k}K(\F_k^i\boxtimes\F_r^j\boxtimes\F_s^l)
\end{eqnarray*}

Putting all these together and expanding the curvature of $\W$ we obtain $K(\W)=0$. \fim

\begin{coro}\label{lema4}
    Let $\W=\F\boxtimes\W_1\boxtimes\cdots\boxtimes\W_n$ be a product of a foliation and $n$ webs. If for each $i,j,l=1,\hdots,n$ the Proposition \ref{prop4} hypotheses are satisfies and we have:
    \begin{enumerate}[i.]
        \item $K(\F\boxtimes\W_i)=0$;
        \item $K(\F\boxtimes\W_i\boxtimes\W_j)=0$, $i\ne j$.
    \end{enumerate}
    Then $K(\W)=0$.
\end{coro}

\dem
Writing as before $\W_i=\F_1^i\boxtimes\cdots\boxtimes\F_{n_i}^i$, we have that
\begin{eqnarray*}
    0 = K(\F\boxtimes\W_i)&=& K(\W_i) + \sum\limits_{1\le k<r\le n_i}K(\F\boxtimes\F_k^i\boxtimes\F_r^i) = \sum\limits_{1\le k<r\le n_i}K(\F\boxtimes\F_k^i\boxtimes\F_r^i)
\end{eqnarray*}
and
\begin{eqnarray*}
    0=K(\F\boxtimes\W_i\boxtimes\W_j)&=&K(\W_i)+K(\W_j) + \sum\limits_{\overset{1\le k < r \le n_i}{1 \le s \le n_j}} K(\F_k^i\boxtimes\F_r^i\boxtimes\F_s^j) \\
    &&+ \sum\limits_{\overset{        1\le k < r \le n_j}{1 \le s \le n_i}} K(\F_k^j\boxtimes\F_r^j\boxtimes\F_s^i) + \sum\limits_{1\le k<r\le n_i}K(\F\boxtimes\F_k^i\boxtimes\F_r^i) \\
    && + \sum\limits_{1\le k<r\le n_i}K(\F\boxtimes\F_k^j\boxtimes\F_r^j) + \sum\limits_{\overset{1\le k<\le n_i}{1\le r \le n_j}} K(\F\boxtimes\F_k^i\boxtimes\F_r^j)\\
    &=& \sum\limits_{\overset{1\le k\le n_i}{1\le r \le n_j}} K(\F\boxtimes\F_k^i\boxtimes\F_r^j)
\end{eqnarray*}
Now a simple computation shows that $K(\W)=0$. \fim

\begin{coro}\label{lema5}
    Let $\W=\F_1\boxtimes\F_2\boxtimes\W_1\boxtimes\cdots\boxtimes\W_n$ be a product of two foliations and $n$ webs. With the same hypotheses as the previous corollary for $\F_1$ and $\F_2$, and additionally: if  $K(\F_1\boxtimes\F_2\boxtimes\W_i)=0$ for each $i=1,\hdots,n$. Then $K(\W)=0$.
\end{coro}

\dem
With the same notation as before, we have

\begin{eqnarray*}
0=K(\F_1\boxtimes\F_2\boxtimes\W_i)&=& K(\W_i) + \sum\limits_{1\le k<r\le n_i} K(\F_1\boxtimes\F_k^i\boxtimes\F_r^i) + \sum\limits_{1\le k<r\le n_i} K(\F_2\boxtimes\F_k^i\boxtimes\F_r^i) \\
&&+ \sum\limits_{k=1}^{n_i} K(\F_1\boxtimes\F_2\boxtimes\F_k^i)\\
&=& \sum\limits_{k=1}^{n_i} K(\F_1\boxtimes\F_2\boxtimes\F_k^i)
\end{eqnarray*}

Combining this with the information obtained in the previous lemmas, we conclude $K(\W)=0$. \fim

Before starting to construct the main results, first note that in the case where we have a foliation and a not invariant line, the dual web will not be flat in general.
 
\subsection{Product with not invariant line}

Given a foliation $\F$, it is expected that for a generic non-invariant line $l$ the web $\leg (l \boxtimes \F)$ is not flat. To simplify the presentation, we will show this in the following case.

\begin{prop}
Let $\F$ be a foliation of degree $d=2$ and $l$ a generic line, then the dual web of $\W=l \boxtimes \F$, $\leg \W$ is not flat.
\end{prop}

The proof will be done by contradiction. Since being flat is a closed condition, we shall assume that $\leg \W$ is flat for all $l \in \Check{\mathbb{P}}^2$.

We can assume that:
\begin{enumerate}[i.]
    \item $l\cap\sing\F=\emptyset$, that is $\Check{\Sigma}^{l}_\F=\emptyset$
    \item $\G_\F(l)\not\subseteq\Check{\Sigma}^{rad}_\F$, in fact, for $p$ a tangency point between $l$ and $\F$ we have that $\G_\F(p)=\Check{l} \notin \Check{\Sigma}^{rad}_\F$.
\end{enumerate}

 Set $\G_\F(l)=C$ and take a generic point $t\in C$ such that the line $\Check{t}$ is not invariant by $\F$. Let $\tang(\F,\Check{t})=p_1+p_2$, with $p_1 \in \Check{t}$, since $I(\F)$ is formed by invariant lines we have that $p_1$ and $p_2$ are not inflection points. On the order hand, condition (ii) guarantees that these points are not radial singularities either. Therefore, in a neighborhood of $t$ we can write
    $$\leg(\W)=(\leg l \boxtimes \F_1)\boxtimes\F_2$$
    where $\F_1$ is tangent to $\leg l$ along $C$ and $\F_2$ is transversal to $\W_2:=\leg l \boxtimes \F_1$

    By \cite[Theorem 1]{JD}, since $K(\leg\W)$ is holomorphic in a generic point of $C$ we have that $C$ is invariant by $\W_2$ or by $\beta_{\W_2}(\F_2)=\F_2$. In the first case $C$ is $\leg l$-invariant and $l$ would be invariant by $\F$, so this is not the case. Thus $C$ is invariant by $\F_2$ and so $\Check{C}$ is an invariant curve for $\F$. Clearly $C$ depends of the line $l$ and then so will $\Check{C}$. Since $\F$ can have at most a pencil of invariant curves we would have an infinite number of lines $l$ with the same image $\G_\F(l)$, which is impossible.

%%%%%%%%%%%%%%%%%%%%%%%%%%%%%%%%%%%%%%%%%%%%%%%%%%%%%%

\section{Convex Reduced Foliations}
In this section we study webs that are dual to the product of convex reduced foliations and (maybe) some lines. We say that $\F$ is \textbf{convex} if $I(\F)$ is formed by invariant lines. If moreover $I(\F)$ is a reduced divisor we will say that $\F$ is \textbf{convex reduced}.

\begin{prop}\label{coro1}
Let $\F_1$, $\F_2$ be convex reduced foliations on $\mathbb{P}_{\mathbb{C}}^{2}$ such that, $\tang(\F_1,\F_2)$ is a union of invariant lines by $\F_1$ and $\F_2$. Then the web $\W=\F_1 \boxtimes \F_2$ satisfies

$$\Delta(\leg\W)=\Check{\Sigma}_{\F_1}^{\text{rad}}\cup \Check{\Sigma}_{\F_2}^{\text{rad}}.$$
\end{prop}

\dem As $\F_1$, $\F_2$ are convex and reduced, we have that all the singularities are non-degenerate (\cite[Lemma 6.8]{B}), thus  $\Delta(\leg\F_i)=\Check{\Sigma}_{\F_i}^{\text{rad}}$ for each $i=1,2$. furthermore, being  $\tang(\F_1,\F_2)=l_1 \cup \hdots \cup l_k$, with each $l_j$ invariant, and due to the fact that the Gauss map of invariant lines being points, we obtain the desired equality as a consequence of Proposition \ref{teo1}. \fim

Now we establish the flatness of the product of two convex reduced foliations under certain hypothesis. 

\begin{prop}\label{teo2}
Let $\F_1$, $\F_2$ be convex reduced foliations on $\mathbb{P}_{\mathbb{C}}^{2}$ of respective degrees $d_1, d_2 \geq 3$, such that, $\tang(\F_1,\F_2)$ in a union of invariant lines and denote $\W=\F_1 \boxtimes \F_2$. Then, the web $\leg\W=\leg\F_1 \boxtimes \leg\F_2$ is flat.
\end{prop}

\dem By Proposition \ref{coro1}, we have that $\Delta(\leg\W)=\Check{\Sigma}_{\F_1}^{\text{rad}}\cup \Check{\Sigma}_{\F_2}^{\text{rad}}$. Take $\Check{s} \in \Check{\Sigma}_{\F_1}^{\text{rad}}\setminus \Check{\Sigma}_{\F_2}^{\text{rad}}$ and a generic line $l$ (particularly, not invariant) passing through $s$. By \cite[Proposition 3.3]{JD} we can write locally around $l$: $\leg\F_1=\W_{\nu_1}\boxtimes\W_{d_1-\nu_1}$ with $l$ a totally invariant component of $\Delta(\W_{\nu_1})$ of minimal multiplicity $\nu_1 -1$ and $\W_{d_1-\nu_1}$ transverse to $\Check{s}$. As $s \not\in \Sigma_{\F_2}^{\text{rad}}$ we have that $\tang(\F_2,l,s)\le 1$.

\begin{itemize}
    \item If $s \not\in \sing\F_2$, we can assume $l \cap \sing\F_2=\emptyset$, then $\leg\F_2$  is completely transversal to $\Check{s}$. On the other hand, if one direction of $\W_{d_1-\nu_1}$ and $\leg\F_2$ coincide, this direction would correspond to a point of  $\tang(\F_1,\F_2)=l_1 \cup \hdots\cup l_k$ therefore, we would have $l=l_j$ for some $j=1,\hdots, k$, which does not occur since $l$ is generic. Thus, $\leg\F_2$ is transverse to $\W_{d-\nu_1}$ and $\leg\W=\W_{\nu_1}\boxtimes(\W_{d_1-\nu_1}\boxtimes\leg\F_2)$ with $\W_{d_1-\nu_1}\boxtimes\leg\F_2$ transverse to $\Check{s}$. Therefore, by \cite[Proposition 2.6]{JD}, the curvature of $\leg\W$ is holomorphic along $\Check{s}$.

    \item If $s \in \sing\F_2$ we have $\tang(\F_2,l,s)=1$ and we decompose $\leg\F_2=\W_1\boxtimes\W_{d_2-1}$. By the previous argument, $\W_{d_2-1}$ is transverse to $\W_{d_1-\nu_1}$. Thus, by \cite[Proposition 3.9]{S}, the curvature of $\leg\W=(\W_1\boxtimes\W_{\nu_1})\boxtimes(\W_{d_1-\nu_1}\boxtimes\W_{d_2-1})$ is holomorphic along $\Check{s}$.
\end{itemize}

Assume now that $s\in \Sigma_{\F_1}^{\text{rad}}\cap \Sigma_{\F_2}^{\text{rad}}$. In this case we can locally write $\leg \W = \W_{\nu_1}\boxtimes \W_{\nu_2}\boxtimes\W_{d_1+d_2-\nu_1 - \nu_2}$ where $\W_{d_1+d_2-\nu_1 - \nu_2}$ is smooth and transverse to $\Check{s}$. We conclude by applying Lemma $ \ref{lema1}$ that $K(\leg \W)$ is holomorphic.
\fim

\begin{ex}\label{ex2}
    Let $\F_d$ be the Fermat foliation of degree $d$, it is defined by the vector field
    $$V_d=(x^d-x)\frac{\partial}{\partial x} + (y^d-y)\frac{\partial}{\partial y}.$$

By \cite[Section 5]{JD} we have that $\F_d$ is convex reduced and that the $3d$ invariant lines are: $L_{\infty},  x=0,  y=0,  x=\xi, y=\xi, y=\xi x$,  with $\xi^{d-1}=1.$

Take now $\F_l$, $\F_{d}$ Fermat foliations of degrees $l$ and $d$ with $d=2l-1$. Note that, $\deg(\tang(\F_l,\F_d))$ $=d+l+1=3l$. On the other hand, with $U_{l-1}(\mathbb{C})$ being the unity group composed of $(l-1)$ roots of unity, we have $U_{l-1}(\mathbb{C})\subseteq U_{2(l-1)}(\mathbb{C})=U_{d-1}(\mathbb{C})$ thus, $\mathcal{I}_{\F_l}\subseteq\mathcal{I}_{\F_d}$. Therefore, $\tang(\F_l,\F_d)=\mathcal{I}_{\F_l}$ contains only invariant lines of both foliations. We conclude by Proposition \ref{teo2} that $\leg (\F_d \boxtimes \F_l)$ is flat.
\end{ex}

The following proposition shows that this is the only way to obtain such examples with Fermat foliations. 

\begin{prop}\label{prop1}
    Let $\F_l$, $\F_d$ be Fermat foliations, with $l<d$. Then, $\tang(\F_l,\F_d)$ is reduced and formed by invariant lines of $\F_l$ and $\F_d$ if and only if $d=2l-1$.
\end{prop}

\dem
If $\tang(\F_l,\F_d)=l_1\cup\hdots\cup l_s$ where $l_j$ are distinct invariant lines of $\F_l$ and $\F_d$, we know that $\F_l$ and $\F_d$ have 
$3l$ and $3d$ invariant lines, respectively, being them
$$L_{\infty}, \hspace{0,3cm} x=0, , \hspace{0,3cm} y=0, \hspace{0,3cm} x=\xi_l, \hspace{0,3cm} y=\xi_l, \hspace{0,3cm} y=\xi_l x, \hspace{1cm} \text{with } \xi_l^{l-1}=1$$
$$L_{\infty}, \hspace{0,3cm} x=0, , \hspace{0,3cm} y=0, \hspace{0,3cm} x=\xi_d, \hspace{0,3cm} y=\xi_d, \hspace{0,3cm} y=\xi_d x, \hspace{1cm} \text{with } \xi_d^{d-1}=1$$
Let $k=gcd(l-1,d-1)$, then $\F_l$ and $\F_d$ have $3k+3$ common invariant lines and we can write $l-1=k\cdot l_1$ and $d-1=k\cdot d_1$ with $gcd(l_1,d_1)=1$. As $\deg(\tang(\F_l,\F_d))=d+l+1$ we have that $d+l=3k+2$, thus $l_1 + d_1=3$. Therefore, $l=k+1$ and $d=2k+1$ i.e $d=2l-1$.\\
From Example \ref{ex2}, we have the other implication. \fim

In a similar way to Proposition \ref{teo2} we can prove the result with three foliations.

\begin{prop}\label{teo7}
    Let $\F_1$, $\F_2$, $\F_3$ be convex reduced foliations on $\mathbb{P}_{\mathbb{C}}^{2}$ of respective degrees $d_1$, $d_2$ and $d_3$, such that, $\tang(\F_i,\F_j)= C_{ij}^{inv}$ where $\forall i,j=1,2,3$, $C_{ij}^{inv}$ is a product of invariant lines.  Consider the web $\W=\F_1\boxtimes\F_2\boxtimes\F_3$. Then $\leg\W=\leg\F_1 \boxtimes \leg\F_2\boxtimes\leg\F_3$ is flat. 
\end{prop}

\dem
As before the discriminant $\Delta(\leg\W)$ is
$$\Delta(\leg\W)=\Check{\Sigma}_{\F_1}^{rad}\cup\Check{\Sigma}_{\F_2}^{rad}\cup\Check{\Sigma}_{\F_3}^{rad}$$

Take $s \in \Sigma_{\F_1}^{rad}\setminus (\Sigma_{\F_2}^{rad}\cup\Sigma_{\F_3}^{rad})$, we have that $s$ is a radial singularity of order $\nu_1-1$ of $\F_1$ and given a generic line $l$ (particularly, not invariant) on $\Check{s}$, we can write $\leg\F_1=\W_{\nu_1}\boxtimes\W_{d_1-\nu_1}$ with $\W_{\nu_1}$ leaving $\Check{s}$ invariant and with $\Delta(\W_{\nu_1})$ having minimal multiplicity along $\Check{s}$ and the $(d_1-\nu_1)$-web $\W_{d_1-\nu_1}$ is transverse to $\Check{s}$ . As $s \not\in \Sigma_{\F_2}^{rad}\cup\Sigma_{\F_3}^{rad}$ we have that $\tang(\F_2,l,s)$ and $\tang(\F_3,l,s)$ are less or equal than $1$.

\begin{itemize}
    \item If $s \not\in (\sing\F_2\cup\sing\F_3)$, we can assume $l \cap (\sing\F_2\cup\sing\F_3)=\emptyset$, then $\leg\F_2$ and $\leg\F_3$ are completely transversal to $\Check{s}$. On the other hand, for example, if one direction of $\W_{d_1-\nu_1}$ and $\leg\F_2$ coincide, this direction would correspond to a point of  $\tang(\F_1,\F_2)= l_1 \cup \hdots\cup l_k$ therefore, we would have $l=l_j$ for some $j$, which does not occur. Thus, $\leg\F_2$ and $\leg\F_3$ are transverse to $\W_{d-\nu_1}$ and in the same way they are mutually transversal, thus in the neighborhood of $l$ we can write $\leg\W=\W_{\nu_1}\boxtimes(\W_{d_1-\nu_1}\boxtimes\leg\F_2\boxtimes\leg\F_3)$ with $\W_{d_1-\nu_1}\boxtimes\leg\F_2\boxtimes\leg\F_3$ smooth and transverse to $\Check{s}$. Therefore, by \cite[Proposition 2.6]{JD}, the curvature of $\leg\W$ is holomorphic along $\Check{s}$.

    \item If $s \in (\sing\F_2\setminus\sing\F_3)$ we have $\tang(\F_2,l,s)=1$ and we decompose $\leg\F_2=\F_2\boxtimes\W_{d_2-1}$ where $\F_2$ is a foliation tangent to $\Check{s}$ and $\W_{d_2-1}$ is a web transverse to $\Check{s}$. By the previous argument, $\W_{d_2-1}$ and $\leg\F_3$ are transverse to $\W_{d_1-\nu_1}$ and to each other as well. Thus, by \cite[Proposition 3.9]{S}, the curvature of $\leg\W=(\F_2\boxtimes\W_{\nu_1})\boxtimes(\W_{d_1-\nu_1}\boxtimes\W_{d_2-1}\boxtimes\leg\F_3)$ is holomorphic along $\Check{s}$.

    \item If $s \in (\sing\F_2\cap\sing\F_3)$ we have $\tang(\F_2,l,s)=\tang(\F_3,l,s)=1$ and we decompose $\leg\F_2=\F_2\boxtimes\W_{d_2-1}$ and $\leg\F_3=\F_3\boxtimes\W_{d_3-1}$ as before. By the previous argument, $\W_{d_2-1}$ and $\W_{d_3-1}$ are transverse to $\W_{d_1-\nu_1}$. Thus, by \cite[Remark 3.10]{S}, the curvature of $\leg\W=(\F_2\boxtimes\F_3\boxtimes\W_{\nu_1})\boxtimes(\W_{d_1-\nu_1}\boxtimes\W_{d_2-1}\boxtimes\W_{d_3-1})$ is holomorphic along $\Check{s}$.
\end{itemize}

Now, take $s \in (\Sigma_{\F_1}^{rad}\cap \Sigma_{\F_2}^{rad})\setminus \Sigma_{\F_3}^{rad}$, we have that $s$ is a radial singularity of orders $\nu_1-1$ and $\nu_2-1$ of $\F_1$ and $\F_2$ respectively and we have local decompositions as before: $\leg \F_i= \W_{\nu_i}\boxtimes\W_{d_i-\nu_i}$. We need to analyze two situations:

\begin{itemize}
    \item If $s \not\in\sing\F_3$ by Lemma \ref{lema1}, we obtain that the curvature of $\leg\W=(\W_{\nu_1}\boxtimes\W_{\nu_2})\boxtimes(\leg\F_3\boxtimes\W_{d_1-\nu_1}\boxtimes\W_{d_2-\nu_2})$ is holomorphic along $\Check{s}$

    \item if $s\in\sing\F_3$, we can decompose $\leg\F_3=\F_3\boxtimes\W_{d_3-1}$ with $\F_3$ a foliation tangent to $\Check{s}$. Again by Lemma \ref{lema2} the curvature of $\leg\W=(\W_{\nu_1}\boxtimes\W_{\nu_2}\boxtimes\F_3)\boxtimes(\W_{d_1-\nu_1}\boxtimes\W_{d_2-\nu_2}\boxtimes\W_{d_3-1})$ is holomorphic along $\Check{s}$.
\end{itemize}

In the case where $\Check{s}$ is radial singularity of $\F_1$, $\F_2$ and $\F_3$ of respective orders $\nu_1-1$, $\nu_2-1$, $\nu_3-1$ we obtain that $\leg\W$ is holomorphic along $\Check{s}$ by Lemma \ref{lema1}. \fim

Now we consider product of two foliations and one invariant line.

\begin{prop}\label{teo3}
Let $\F_1$, $\F_2$ be convex reduced foliations on $\mathbb{P}_{\mathbb{C}}^{2}$ of respective degrees $d_1$ and $d_2$, such that, $\tang(\F_1,\F_2)$ is a union of invariant lines. For any invariant line $l$ consider the web $\W=l\boxtimes\F_1\boxtimes\F_2$. Then $\leg\W=\leg l\boxtimes\leg\F_1 \boxtimes \leg\F_2$ is flat.
\end{prop}

\dem
It follows from Proposition \ref{coro1} and \cite[Section 2]{S} that $\Delta(\leg\W)=\Check{\Sigma}_{\F_1}^{\text{rad}}\cup \Check{\Sigma}_{\F_2}^{\text{rad}}\cup\Check{\Sigma}_{\F_1}^{l}\cup\Check{\Sigma}_{\F_2}^{l}$, where $\Sigma^l_{\F}= \sing \F \cap l$. Consider the following cases

i. $s\in\Sigma_{\F_1}^{\text{rad}}\setminus(\Sigma_{\F_2}^{\text{rad}}\cup l)$ then we repeat the argument of Proposition \ref{teo2} and obtain that $K(\leg \W)$ is holomorphic along $\Check{s}$. The case in which $\Check{s}\in\Check{\Sigma}_{\F_2}^{\text{rad}}\setminus(\Check{\Sigma}_{\F_1}^{\text{rad}}\cup l)$ is analogous.\\

ii. If $s\in(\Sigma_{\F_1}^{\text{rad}}\cap\Sigma_{\F_2}^{\text{rad}})\setminus l$ with radiality orders $\nu_1-1$ and $\nu_2-1$ respectively. Since $s \not\in l$ we obtain that the direction of $\leg l$ on the general point of $\Check{s}$ does not coincide with any other of $\leg\F_1$ and $\leg\F_2$. So we can write 
$$\leg\W=(\W_{\nu_1}\boxtimes\W_{\nu_2})\boxtimes(\leg l\boxtimes\W_{d_1+d_2-\nu_1-\nu_2})$$
where $\leg l\boxtimes\W_{d_1+d_2-\nu_1-\nu_2}$ is smooth and transverse to $\Check{s}$. Again by Lemma \ref{lema1} we have that $K(\leg\W)$ is holomorphic on $\Check{s}$.\\

iii. If $s\in l\setminus\Sigma_{\F_2}^{\text{rad}}$ with $\nu_1-1$ being the radiality order of $s$ in $\F_1$, then we write $\leg\F_1=\W_{\nu_1}\boxtimes\W_{d_1-\nu_1}$. Note that the directions of $\leg\F_2$ and $\W_{d_1-\nu_1}$ do not coincide at generic $\Check{t}\in \Check{s}$, and $\tang(\F_2,t,s)\le 1$ (because $s\not\in\Sigma_{\F_2}^{\text{rad}}$). Then,
\begin{itemize}
    \item If $s \notin \sing\F_2$, we have that,
    $$\leg\W=(\leg l\boxtimes\W_{\nu_1})\boxtimes(\leg\F_2\boxtimes\W_{d_1-\nu_1})$$
    where $\leg\F_2\boxtimes\W_{d_1-\nu_1}$ is smooth and transversal to $\Check{s}$. Now, we apply  \cite[Proposition 3.9]{S} if $\nu_1 \ge 2$ and \cite[Theorem 1]{JD} if $\nu_1=1$ in order to obtain the flatness.
    
    \item If $s \in \sing\F_2$, we have $\tang(\F_2,t,s)=1$ and we can write
    $$\leg\W=(\leg l\boxtimes\W_1\boxtimes\W_{\nu_1})\boxtimes(\W_{d_1-\nu_1}\boxtimes\W_{d_2-1})$$
    with $\W_{d_1-\nu_1}\boxtimes\W_{d_2-1}$ smooth and transverse to $\Check{s}$ and $\W_1$ a foliation tangent to $\Check{s}$. If $\nu_1\ge 2$ the flatness of $\leg \W$ will follow from \cite[Remark 3.10]{S}. On other hand, if $\nu_1=1$ we conclude by \cite[Theorem 1.1]{M}.
\end{itemize}

iv. If $s\in\Sigma_{\F_1}^{\text{rad}}\cap\Sigma_{\F_2}^{\text{rad}}\cap l$ with $\nu_1-1$ and $\nu_2-1$ being the orders of $s$, respectively, in $\F_1$ and $\F_2$, then we can write
$$\leg\W=(\leg l\boxtimes\W_{\nu_1}\boxtimes\W_{\nu_2})\boxtimes(\W_{d_1+d_2-\nu_1-\nu_2})$$
with $\W_{d_1+d_2-\nu_1-\nu_2}$ smooth and transverse to $\Check{s}$. We conclude by Lemma \ref{lema2}. \fim

\noindent\underline{Proof. of Theorem A}:

Set $\W'=\leg l_i \boxtimes\cdots\boxtimes \leg l_k$ and $\W''=\leg\F_1\boxtimes\cdots\boxtimes\leg\F_n$. By \cite[lemma 2.1]{S2}, we have that
{\small$$K(\W'\boxtimes\W'')=K(\W') -(k-2)\sum_{i=1}^k K(\leg l_i\boxtimes\W'')+\sum_{1\le i<j\le k} K(\leg l_i\boxtimes\leg l_j\boxtimes\W'') + \left ( \begin{matrix}
    k-1 \\ 2
\end{matrix} \right )K(\W'').$$} 

By \cite[Theorem 4.2]{JD}, Propositions \ref{teo2}, \ref{teo7}, and \ref{teo3} we are in condition to apply Corollary \ref{lema4} and conclude that $K(\leg l_i\boxtimes\W'')=0$. Similarly, we use in addition \cite[Theorem 2]{S2} to be able to use Corollary \ref{lema5} and obtain $K(\leg l_i\boxtimes\leg l_j\boxtimes\W'')=0$. Finally, by Proposition \ref{prop4} we have $K(\W'')=0$ and we finish our proof.
\fim

\begin{ex}
    Take $\F_d$ and $\F_{2d-1}$ Fermat foliations of degree $d$ and $2d-1$. By Proposition \ref{prop1}, the curve $\tang(\F_d,\F_{2d-1})$ is formed by the invariant lines: 
    $$L_{\infty}, \hspace{0,3cm} x=0, , \hspace{0,3cm} y=0, \hspace{0,3cm} x=\xi, \hspace{0,3cm} y=\xi, \hspace{0,3cm} y=\xi x, \hspace{1cm} \text{with } \xi^{d-1}=1.$$
    
    In this way, if $C$ is the reduced curve formed by the product of some of these lines and $\W=C\boxtimes\F_d\boxtimes\F_{2d-1}$, we have that $\leg(\W)$ is flat.
\end{ex}

%%%%%%%%%%%%%%%%%%%%%%%%%%%%%%%%%%%%%%%%%%%%%%%%%%%%

\section{Homogeneous Foliations}
This section is devoted to prove Theorem B. We begin by a extension of \cite[Lemma 3.1]{M} for product of non-saturated homogeneous foliations.

\begin{lema}\label{lema3}
    Let $\W=l_1\boxtimes\cdots\boxtimes l_k\boxtimes\h_1 \boxtimes \cdots \boxtimes\h_n$ be the product of homogeneous foliations and lines passing through the origin. If the curvature of $\leg\W$ is holomorphic on $\Check{\mathbb{P}}^2\setminus \Check{O}$, then $\leg\W$ is flat. 
\end{lema}

\begin{dem}
    Fix coordinates $(a,b)$ in $\Check{\mathbb{P}}^2$ associated to the line $\{ax+by=1\}$ in $\mathbb{P}^2$. Set $\omega_i=A_idx+B_idy$ the homogeneous form defining $\h_i$ and $l_j=\{a_jx+b_jy=0\}$.
    
    Thus the web $\W$ is defined by the form $\Omega=F_k\cdot\omega_1\cdots\omega_n$, with $F_k(x,y)=(a_1x + b_1y) \cdot \ldots \cdot(a_kx +b_k y)$. Note that $\W$ is invariant by the homothecies $h_\lambda(x,y)=(\lambda x,\lambda y)$, in fact
    $$h_\lambda^{*}(\Omega)=(\lambda^kF_k)\cdot(\lambda^{d_1+1}\omega_1)\cdots(\lambda^{d_n+1}\omega_n)=\lambda^{d+n+k}\Omega$$
    where $d=d_1+\hdots+d_n$. So we have that $\leg\W$ is invariant by the dual maps $\Check{h}_\lambda(a,b)=(\frac{a}{\lambda},\frac{b}{\lambda})$ and so
    $$\Check{h}_\lambda^{*}(K(\leg\W))=K(\leg\W).$$
    How in affine coordinates $(a,b)$ the $K(\leg\W)$ is a rational $2$-form with poles at maximum in $\Check{O}=L_{\infty}$, we can write the curvature as $K(\leg\W)=p(a,b)da\wedge db$ where $p$ is a polynomial. Therefore $\lambda^2p(a,b)=p(\frac{a}{\lambda},\frac{b}{\lambda})$ and the lemma follows. \fim
\end{dem}

\begin{prop}\label{prop3.2}
    Let $\h_1$, $\h_2$ be convex homogeneous foliations on $\mathbb{P}^{2}$ such that, $\tang(\h_1,\h_2)=L_\infty \cup l_1 \cup \hdots \cup l_k$ where $l_j$ is an invariant line $\forall j=1,\hdots,k$. Consider the web $\W=\h_1\boxtimes\h_2$. Then $\leg\W=\leg\h_1 \boxtimes \leg\h_2$ is flat. 
\end{prop}

\dem
It follows from Proposition \ref{teo1} and \cite[Lemma 3.2]{B} that 
$$\Delta(\leg\W)=\Check{\Sigma}_{\h_1}^{rad}\cup\Check{\Sigma}_{\h_2}^{rad}\cup\Check{O}$$
If $s\in(\Sigma_{\h_1}^{rad}\cup\Sigma_{\h_2}^{rad})$, we do a similar analysis to that done in the Proposition \ref{teo2} and obtain that $K(\leg\W)$ is holomorphic over $\Check{s}$. Thus $K(\leg \W)$ is holomorphic on $\Check{\mathbb{P}}^2_\mathbb{C}\setminus\Check{O}$ and by Lemma \ref{lema3} we have that $\leg\W$ is flat. \fim

We can prove a similar result with three foliations.

\begin{prop}\label{teo6}
    Let $\h_1$, $\h_2$, $\h_3$ be convex homogeneous foliations on $\mathbb{P}_{\mathbb{C}}^{2}$ of respective degrees $d_1$, $d_2$ and $d_3$, such that, $\tang(\h_i,\h_j)=L_\infty \cup C_{ij}^{inv}$ where $\forall i,j=1,2,3$, $C_{ij}^{inv}$ is a product of invariant lines.  Consider the web $\W=\h_1\boxtimes\h_2\boxtimes\h_3$. Then $\leg\W=\leg\h_1 \boxtimes \leg\h_2\boxtimes\leg\h_3$ is flat. 
\end{prop}

\dem
As before the discriminant $\Delta(\leg\W)$ is
$$\Delta(\leg\W)=\Check{\Sigma}_{\h_1}^{rad}\cup\Check{\Sigma}_{\h_2}^{rad}\cup\Check{\Sigma}_{\h_3}^{rad}\cup\Check{O}$$
If $s \in (\Sigma_{\h_1}^{rad}\cup \Sigma_{\h_2}^{rad}\cup\Sigma_{\h_3}^{rad})$, we can proceed as in the Proposition \ref{teo7} and obtain that $K(\leg\W)$ is holomorphic over $\Check{s}$. So since $K(\leg\W)$ is holomorphic on $\mathbb{P}^2_{\mathbb{C}}\setminus\Check{O}$ by Lemma \ref{lema3} we  conclude that $\leg\W$ is flat. \fim

Now we can also consider the product with invariant lines. 

\begin{prop}\label{teo4}
Let $\h_1$, $\h_2$ be convex homogeneous foliations on $\mathbb{P}_{\mathbb{C}}^{2}$ of respective degrees $d_1$ and $d_2$, such that, $\tang(\h_1,\h_2)=L_\infty \cup l_1 \cup \hdots \cup l_k$ where  $l_j$ is a line invariant $\forall j=1,\hdots,k$, and let $l$ be an invariant line. Consider the web $\W=l\boxtimes\h_1\boxtimes\h_2$. Then $\leg\W=\leg l\boxtimes \leg\h_1 \boxtimes \leg\h_2$ is flat. 
\end{prop}

\dem
Firstly, we consider $l=L_\infty$. By \cite{B}, we can locally decompose $\leg(L_\infty\boxtimes\h_1\boxtimes\h_2)$ as
$$\leg(L_\infty\boxtimes\h_1\boxtimes\h_2)=\leg(L_\infty)\boxtimes\W_{d_1+d_2}$$
where $\W_{d_1+d_2}=\F_1\boxtimes\cdots\boxtimes\F_{d_1+d_2}$ and, for each $i$, $\F_i$ is given by $\Check{\omega}_i:=\lambda_i(p)dq-qdp$ with $\lambda_i(p)=p-p_i(p)$ and $\{p_i(p)\}=\underline{\G}_\h^{-1}(p)$. The radial foliation $\leg(L_\infty)$ is defined by the vector field $X:=q\frac{\partial}{\partial q}$ and this is a transverse symmetry of the web $\W_{d_1+d_2}$. By \cite[Theorem 3.4]{S} and Proposition \ref{prop3.2}, we have that $K(\leg(L_\infty\boxtimes\h_1\boxtimes\h_2))=K(\leg(\h_1\boxtimes\h_2))=0$.\\

Now, assume $l$ is an invariant line passing through the origin, we can decompose the discriminant  
{\small $$\Delta(\leg\W)=\Delta(\leg\h_1) \cup \Delta(\leg\h_2)\cup\tang(\leg l,\leg\h_1)\cup\tang(\leg l,\leg \h_2)\cup\tang(\leg\h_1,\leg\h_2)$$}

Since $l$ is invariant by $\h_1$ and $\h_2$, then $\tang(\leg l, \leg\h_1)=\tang(\leg l,\leg\h_2)=\Check{O}\cup\Check{a}$ with $a=l\cap L_\infty$. On the other hand, since $\tang(\h_1,\h_2)$ is a union of invariant lines, we have that 
$$\Delta(\leg\W)=\Check{\Sigma}_{\h_1}^{rad}\cup\Check{\Sigma}_{\h_2}^{rad}\cup\Check{O}\cup\Check{a}.$$
Note that, for $s\in(\Sigma_{\h_1}^{rad}\cup\Sigma_{\h_2}^{rad})\setminus l$ a similar analysis to that done in Proposition \ref{teo3} leads us the holomorphy of $K(\leg \W)$ along $\Check{s}$. Take now a generic line $t$ (particularly, $t$ is not invariant by $\h_1$ and $\h_2$) passing through $a$, we have three situations:
\begin{itemize}
    \item If $a\not\in(\Sigma_{\h_1}^{rad}\cup\Sigma_{\h_2}^{rad})$, since $\h_1$ and $\h_2$ are convex homogeneous, $\tang(\h_1,t,a)=\tang(\h_2,t,a)=1$ thus we can write locally around $\Check{t}$: $\leg\h_1=\F_1\boxtimes\W_{d_1-1}$ and $\leg\h_2=\F_2\boxtimes\W_{d_2-1}$ with $\Check{a}$ invariant by $\F_1$ and $\F_2$ and $\W_{d_1-1}\boxtimes\W_{d_2-1}$ smooth and transverse to $\Check{a}$ . Therefore
    $$\leg\W=(\leg l\boxtimes\F_1\boxtimes\F_2)\boxtimes(\W_{d_1-1}\boxtimes\W_{d_2-1}).$$
    By \cite[Theorem 1.1]{M}, we conclude that $K(\leg\W)$ is holomorphic along $\Check{a}$.

    \item If $a\in\Sigma_{\h_1}^{rad}\setminus\Sigma_{\h_2}^{rad}$, being $\nu_1$ the radiality order of $a$ as singularity of $\h_1$, we can write locally: $\leg\h_1=\W_{\nu_1}\boxtimes\W_{d_1-\nu_1}$ and $\leg\h_2=\F\boxtimes\W_{d_2-1}$ as before. Thus
    $$\leg\W=(\leg l\boxtimes\W_{\nu_1}\boxtimes\F)\boxtimes(\W_{d_1-\nu_1}\boxtimes\W_{d_2-1})$$
     with $\F\boxtimes\W_{d_1-\nu_1}$ tangent to $\Check{a}$ and $\W_{d_1-\nu_1}\boxtimes\W_{d_2-1}$ transverse to $\Check{a}$. By \cite[Remark 3.10]{S} we get the holomorphy of $K(\leg\W)$ along $\Check{a}$.

    \item The case that $a\in (\Check{\Sigma}_{\h_1}^{rad}\cap\Check{\Sigma}_{\h_2}^{rad})$ is similar to what we did in item i.v. of Proposition \ref{teo3}. 
\end{itemize}

Therefore, $K(\leg\W)$ is holomorphic on $\Check{\mathbb{P}}^2\setminus\Check{O}$ and, $\leg\W$ is flat by Lemma \ref{lema3}.  \fim

\noindent\underline{Proof. of Theorem B}:
With the results of this section added to those presented in the previous section, similarly to what was done in the demonstration of Theorem A, Theorem B is established.\\

\begin{ex}
Let $\omega_d$ be the form $\omega_d=y^d dx-x^d dy$. It is shown in \cite[Proposition 4.1]{B} that this form defines a convex homogeneous foliation of degree $d$ with invariant lines: $L_\infty$, $x=0$, $y=0$, $y=\xi x$, with $\xi^{d-1}=1$.

Take now $\h_1$, $\h_2$, $\h_3$ the foliations defined by $\omega_d$, $\omega_{d+1}$ and $\omega_{d+2}$ respectively, with $d=2k+1 \geq 3$. We have that the tangencies between each pair of foliations are: $\tang(\h_1,\h_2)=x^dy^d(y-x)$, $\tang(\h_1,\h_3)=x^dy^d(y^2-x^2)$, $\tang(\h_2,\h_3)=x^{d+1}y^{d+1}(y-x)$. Therefore, we are in the hypothesis of Theorem B and for any curve $C$ formed by the reduced product of any combination of the lines $L_\infty$, $x=0$, $y=0$ and $y=x$ we conclude that $\leg(C\boxtimes \h_1\boxtimes\h_2\boxtimes\h_3)$ is a flat web.\\
\end{ex}

\begin{ex}\label{ex3}
Let $\F_1$ and $\F_2$ be foliations given by the respective $1$-forms $\omega_1=\lambda ydx + xdy$ and $\omega_2=y^2dx+x^2dy$. Consider the $3$-web $\W=\F_1\boxtimes\F_2$ and note that $\tang(\F_1,\F_2)=L_\infty\cup l_1\cup l_2\cup T$ where the lines $l_1:\{x=0\}$ and $l_2:\{y=0\}$ are invariant by $\F_1$ and $\F_2$, and $T:\{y+\lambda x\}$ is transversal. Take $p=\G_{F_1}(T)$ then clearly $\Check{p}\subset \Delta(\leg\W)$ and it is not difficult to verify, by barycenter criterion (see \cite[Theorem 1.1]{M}), that the curvature of $\leg\W$ is not holomorphic in $\Check{p}$.
\end{ex}

\bibliography{main}

\end{document}